\documentclass{article}

\usepackage{arxiv}

\usepackage[utf8]{inputenc} 
\usepackage[T1]{fontenc}    
\usepackage{hyperref}       
\usepackage{url}            
\usepackage{booktabs}       
\usepackage{amsfonts}       
\usepackage{nicefrac}       
\usepackage{microtype}      
\usepackage{lipsum}		
\usepackage{graphicx}
\usepackage{natbib}
\usepackage{dsfont}
\usepackage{doi}
\usepackage{latexsym}
\usepackage{graphicx}
\usepackage{mathptmx}
\usepackage[T1]{fontenc}
\usepackage{subcaption}
\usepackage{listings}
\usepackage{xcolor}
\usepackage{makecell}
\usepackage{pifont}
\usepackage{amsmath}
\usepackage{amsfonts}
\usepackage{amssymb}
\usepackage{amsbsy}
\usepackage{amsthm}
\usepackage{appendix}
\newtheorem{theorem}{Theorem}[section]

\newtheorem{remark}[theorem]{Remark}

\title{Resilience Quantification and its Support for Operational Resilience}

\author{
    {\hspace{1mm}Ion Matei} \\
	Fujitsu Research of America\\
	\texttt{imatei@fujitsu.com} \\
	\And  
    {\hspace{1mm}Maksym Zhenirovskyy} \\
	Fujitsu Research of America\\
	\texttt{mzhenirovskyy@fujitsu.com} \\
}

\hypersetup{
pdftitle={A template for the arxiv style},
pdfsubject={q-bio.NC, q-bio.QM},
pdfauthor={David S.~Hippocampus, Elias D.~Striatum},
pdfkeywords={First keyword, Second keyword, More},
}

\begin{document}
\maketitle

\begin{abstract}
We present a method to quantify a system’s \emph{resilience capacity}, i.e., the set of degradation magnitudes for which all functional requirements remain satisfied.  These requirements come from human stakeholders (e.g., operators, planners) who define the acceptable performance envelope.  By representing the resilience capacity in degradation space, we obtain an application‑agnostic resilience metric (e.g., capacity volume).  To approximate the capacity efficiently in high‑dimensional spaces, we pair machine‑learning classifiers with entropy‑based active sampling, reducing costly feasibility tests.  The learned model then drives \emph{diagnosis} (current health estimation) and \emph{prognostics} (health‑state forecasting) that estimates useful life. These two steps can be complemented by a \textit{reconfiguration} step implemented by human operators to prolong the system's functionality.  An illustrative case study, i.e., a manufacturing  production line meeting weekly human set part demand, demonstrates the proposed workflow.
\end{abstract}


\section{Introduction}
This article develops a methodology for \emph{quantifying} resilience and using it to predict and, when possible, extend a system’s remaining useful life (RUL) under degradations. Building on power‑system concepts \cite{7842605}, we define resilience as the ability to withstand, adapt to, and recover from disruptive events while continuing to meet functional requirements. The central construct is the \emph{resilience capacity}: the set of degradation vectors from which the system can recover without loss of functionality. In degradation space, this capacity marks a “human‑acceptability” region, as it safeguards requirements like safety limits or throughput quotas defined by human stakeholders. Explicitly tying capacity to user‑defined requirements is critical for human‑machine teaming, ensuring that automated decisions remain transparent and trustworthy. The volume of this capacity induces a quantitative resilience metric. To operationalize this metric, we introduce a \emph{health‑state} vector, parameterizing how degradations (e.g., tool wear) affect system behavior. As long as the health state remains within the learned capacity, functionality is preserved; crossing its boundary signals imminent loss of functionality. Accurate health‑state estimation and forecasting enable proactive measures, such as adjusting task requirements or scheduling maintenance. Together with a decision‑support dashboard, this approach allows reconfiguration recommendations to be implemented by human operators, illustrating how human intervention can augment system resilience while balancing cost, risk, and service quality.

The workflow has offline and online phases. Offline, an \emph{active‑sampling} routine explores the degradation space within a fixed test‑budget, invoking a requirement‑satisfaction (feasibility) oracle to label points as feasible or infeasible. These labeled data train a classifier that approximates the capacity boundary. Online, sensor data feed a \emph{diagnosis} module that updates the health state, while a \emph{prognostics} model forecasts when the degradation trajectory will cross this boundary, yielding RUL estimates. If a violation appears imminent, a \emph{reconfiguration} optimizer can propose schedule adjustments or mode switches to keep the system within its resilience capacity as long as possible. The online reconfiguration step is not explicitly addressed in this paper. More details on the implementation of the operational online phase can be found in \cite{10.1145/3631612}.

\textbf{Contributions and paper structure:} We introduce the \emph{resilience capacity}, a formal volume‑based resilience metric tied to functional requirements. We employ \emph{active learning}, an entropy‑guided sampler augmented with monotonic logic labels, to efficiently trace the high‑dimensional capacity using few oracle queries. We estimate \emph{RUL under uncertainty} by propagating the resilience‑capacity model uncertainty into the RUL prediction. Section \ref{sec:background} outlines the problem and workflow, Section \ref{sec:manufacturing-line} describes the manufacturing case study, Section \ref{sec:res-cap} details the active‑learning approach, Section \ref{sec:phm} integrates diagnosis, prognostics, and Section \ref{sec:discussion} places the approach within the broader context of system resilience.
\section{Problem description}
\label{sec:background}
In this section, we define the resilience capacity and introduce an implicit representation that can be queried. We consider a system $\mathcal S$ that executes function‑related tasks while meeting a set of Boolean requirements $\mathcal R=\{r_1,r_2,\dots r_{n_R}\}$, where $n_R$ is the number of requirements for the system.  
These requirements formalize objectives specified by human stakeholders: plant managers that translate customer
orders into production‑rate bounds, safety engineers that codify regulatory limits, and operators that
impose workload constraints.  A system is deemed \emph{functional} when every requirement $r_i$ evaluates to
\textit{true}.  Degradations are represented by the normalized vector
$d=[d_1,d_2,\ldots, d_{n_d}]\in[0,1]^{n_d}$ (with $n_d$ denoting the number of degradations that can affect the system, e.g., machine production capacity reductions), and the resulting degradation space is the hypercube
$\mathcal D=[0,1]^{n_d}$.  We define the \emph{resilience capacity} \(
  \mathcal C \;=\;\bigl\{ d\in\mathcal D 
          \;\bigl|\; \bigwedge_{r\in\mathcal R} r=\text{true}\bigr\}
\), namely the set of all degradation magnitudes under which \textit{every} human‑defined
requirement remains satisfied.  Because the requirements encode what users perceive as
mission‑critical, $\mathcal C$ corresponds to a “human‑acceptability region’’ in the
degradation space.

Although the geometry of $\mathcal C$ may be convex in some applications, it is generally
difficult to characterize analytically, especially in high dimensions.  We therefore view the
partition $\mathcal D=\mathcal C\cup(\mathcal D\setminus\mathcal C)$ as a binary
classification problem and introduce a classifier
$c:\mathcal D\!\to\!\{0,1\}$ given by $c(d)=1$ for $d\in\mathcal C$, and $c(d)=0$ for $d\in\mathcal D\setminus\mathcal C$. Labels are obtained by querying an expensive “oracle’’ (based on simulation or optimization) that
checks all requirements.  To minimize these costly queries, we exploit the following degradation monotonicity properties that hold for our use case:
\begin{eqnarray}
    \label{eq:03251431}
    {d}\in \mathcal{C} &\implies& \{ \hat{{d}}\in \mathcal{D} | \hat{d}_i\leq d_i \} \subseteq \mathcal{C},\\
    \label{eq:03251433}
    {d}\in \mathcal{D}\setminus \mathcal{C} &\implies& \{ \hat{{d}} \in \mathcal{D} | \hat{d}_i > d_i \}\subseteq \mathcal{D}\setminus \mathcal{C},    
\end{eqnarray}
where indices indicate vector entries. The property (\ref{eq:03251431}) states that if the system is functional under a degradation vector then it is functional for degradations that are less intense. 
These properties can be used to generate additional labels, without the need for querying the oracle.
Thus one oracle call can label infinitely many neighbors, guiding an
entropy‑based sampler toward the most informative points near the boundary of
$\mathcal C$.

After learning a surrogate for $\mathcal C$, we can approximate its volume, which serves as a
scalar resilience metric: a volume close to~1 (the volume of $\mathcal D$) indicates that most
degradation scenarios are tolerable.  While exact volume computation is \#P‑hard (\cite{doi:10.1137/0217060}, Theorem 1), random‑walk
algorithms \cite{10.1145/102782.102783,dadush:LIPIcs.SOCG.2015.704} allow polynomial‑time estimates; here, the learned
classifier plays that oracle role.  Boundary analysis of $\mathcal C$ further identifies the
directions in degradation space that most threaten functionality, informing the prognostics and health management (PHM)
procedures discussed in the next sections.

\section{Manufacturing Line Use Case}
\label{sec:manufacturing-line}
Inspired by a real industrial scenario, we present a use case for demonstrating the proposed methodology, where multiple part types are processed sequentially through three stages, denoted $A$, $B$, and $C$. Each process is supported by a distinct group of specialized machine tools. Machines in stages $A$ and $C$ are restricted to one part type at a time, while machines in stage $B$ may serve multiple types simultaneously. The workflow begins in process $A$, where raw materials are handled and then set aside for an aging period. Subsequent processing stages proceed without delay. Each machine exhibits part-specific production rates. The manufacturing line contains approximately 50 machines of types $A$ and $C$, around 20 machines of type $B$, and supports 40 distinct part types. The machines can be manually, (i.e., by human operators) configured to process specific part types.  Weekly demand forecasts are issued by human planners and must be met through scheduling. These human-issued requirements represent explicit performance targets and embed a human-in-the-loop dimension into the scheduling process. Each machine must be manually configured to process a given part type, and frequent reconfigurations are discouraged due to setup time. Inspired by queuing theory and process-flow models~\cite{10.5555/2693068.2693082}, we formulate the scheduling problem mathematically.

Let $\mu_{i,j}$ be the production rate of machine $i$ when processing part $j$, and let $\theta_i$ indicate the machine’s current mode (including an idle mode), i.e., a physical configuration of machine $i$, enabling it to process a part $j$. Index $i$ ranges from 1 to the total number of machines ($\sim$ 110). Index $j$ ranges from 1 to the number of distinct parts ($\sim$ 40).  Each machine $i$ maintains part-specific queues, whose evolution is governed by $Q_{i,j}[k+1] = Q_{i,j}[k] + \lambda_{i,j}[k] - \bar{\mu}_{i,j}[k]$, where $\lambda_{i,j}$ denotes the part type $j$ input rate at machine $i$, $\bar{\mu}_{i,j}[k] = \mathds{1}_{\{\theta_i[k]=j\}}\,\mu_{i,j}$, with $\mathds{1}$ denoting the indicator function, and $k$ representing the time index. Changing machine $i$'s mode to process type $j$ parts incurs a setup time $T_{i,j}$. Assuming that the control interval $T$ exceeds $T_{i,j}$, the production rate is temporarily reduced after switching: $\mu_{i,j}^s =\mu_{i,j}(T - T_{i,j})/{T}$. Thus, the output rate expression becomes:
\(
\bar{\mu}_{i,j}[k] = \mathds{1}_{\{\theta_i[k]=j,\theta_i[k-1]=j\}}\,\mu_{i,j} + 
\mathds{1}_{\{\theta_i[k]=j,\theta_i[k-1]\neq j\}}\,\mu_{i,j}^s.
\)
The effect of aging in process $A$ is modeled by deadline shifting. For instance, producing 100 parts in 5 days with 1 day aging implies a production window of 4 days.

The process flow model includes three main constraints. The \textit{queue bounds}, i.e., $0 \le Q_{i,j}[k] \le Q_{i,j}^{\max}$, with $Q_{i,j}^{\max}$ a positive scalar, limit the number of accumulated parts of type $j$ at machine $i$. The \textit{input rate bounds}, i.e., $0 \leq \lambda_{i,j}[k] \leq \lambda_{i,j}^{max}$ enforce the non-negativity of the input rates, where the positive scaler $\lambda_{i,j}^{max}$ imposes limits on the available raw materials.   The \textit{flow conservation}, i.e., $\sum_i \bar{\mu}_{i,j}[k] = \sum_l \lambda_{l,j}[k]$, $\forall j$ mathematically states that only parts that have been processed at a previous stage can move to the next stage, and that their total number must be preserved. Here $\bar{\mu}_{i,j}[k]$ represents part $j$ output rate from the machine $i$ at a previous stage (e.g., type $A$ machine), and $\lambda_{i,j}[k]$ denotes part $j$ input rate to machine $l$ in the next stage (e.g., type $B$ machine). Machines in process $B$ may operate on multiple parts in parallel; for modeling simplicity, we treat them as parallel virtual machines, each processing a single type. Let $T_d = n_d T$ denote the deadline for delivering $N_j$ units of part $j$. This demand corresponds to an average output rate of:
\(
\frac{1}{n_d} \sum_{k=1}^{n_d} \bar{\mu}_j^C[k] = \mu_j^d\), with  \(\mu_j^d = N_j/{T_d}
\), where $\bar{\mu}^C_j$ denotes the output rate for part type $j$, produced by type $C$ machines. To check whether this target is achievable under degradation, we solve a mixed-integer linear program (MILP). Let $x_{i,j}[k] = \mathds{1}_{\{\theta_i[k]=j\}}$ be binary variables encoding machine-mode configurations, $\lambda_{i,j}[k]$ the input rates, and $\bar{\mu}_j^C[k]$ the outputs from process $C$. Let bold letters denote the optimization variables over the time horizon $T_d$. We group these variables into a single decision vector $\boldsymbol{z} = (\boldsymbol{Q}, \boldsymbol{x}, \boldsymbol{\lambda}, \boldsymbol{\mu}^C)$ and write:
\begin{eqnarray}
\label{eq:040200920}
\min_{\boldsymbol{z}} \quad & \boldsymbol{c}^T \boldsymbol{z}, \\
\nonumber
\text{s.t.} \quad & \boldsymbol{A}_{eq}\,\boldsymbol{z} = \boldsymbol{b}_{eq}, \\
\nonumber
            \quad & \boldsymbol{A}_{ineq}\,\boldsymbol{z} \leq \boldsymbol{b}_{ineq}, \\
\nonumber
& \boldsymbol{z}_{lb} \le \boldsymbol{z} \le \boldsymbol{z}_{ub},                        
\end{eqnarray}
with the binary variables $x_{i,j} \in \{0,1\}$. The goal of this formulation is to generate a schedule that ensures a target production rate (implemented explicitly as a inequality constraint the enforces the output rates of machines C to be at least as large as a target rate), while minimizing a linear cost.  The cost vector $\boldsymbol{c}$  include manufacturing costs (e.g., energy), and it can be made zero in feasibility checks. The linear equality constraints, defined by matrix $\boldsymbol{A}_{eq}$ and vector $\boldsymbol{b}_{eq}$, enforce queue dynamics and flow conservation. The inequality bounds constrain queue lengths and enforce nonnegativity, while the linear inequalities, defined by matrix $\boldsymbol{A}_{ineq}$ and vector $\boldsymbol{b}_{ineq}$, are used to set the minimal production targets. 
 Given our manufacturing line topology and machine capabilities, each time step introduces approximately 10,000 optimization variables, with roughly half of them being binary variables. The total variable count scales linearly with the horizon length. When the integrality of the machine modes is enforced, the resulting MILP becomes challenging to solve for long time horizons. Then, we either apply variable relaxations or split the time horizon into smaller intervals and solve for each interval, while updating the queues  after each interval solution. Both approaches were used to generate (approximate) schedules.

\section{Resilience Capacity}
\label{sec:res-cap}
 We model machine degradations as a reduction in the maximum achievable production rates. For a machine $i$ processing part $j$, the degraded rate is given by $(1-d_{i})\mu_{i,j}$, where $d_{i} \in [0, 1]$ quantifies the degradation. A value of $d_{i}=0$ indicates full functionality, while $d_{i}=1$ indicates complete failure. Here we assume that the changes in degradations are much slower than machine scheduling interval. The parameters of problem (\ref{eq:040200920}) (e.g., matrices in the linear equalities and inequalities) become functions of the degradation vector $d$. To test if a target production rate is feasible under a given degradation vector, we relax the MILP to an LP  by allowing binary variables to vary in $[0,1]$ and model one aggregated weekly interval. This produces a necessary condition for feasibility. The intuition of this condition is as follows: we need at least one feasible schedule to find (average) production rates.  If however no feasible production rate exists due to the infeasibility of the relaxed LP, then no valid schedule exists. Even this LP has nearly 5,000 variables and thousands of constraints. Solving this LP is equivalent to querying an \textit{oracle} to check feasibility. Solving this problem for various ${d}$ defines the feasible degradation set, or \emph{resilience capacity} $\mathcal{C}$. Under the LP formulation in (\ref{eq:040200920}), it can be shown that $\mathcal{C}$ forms a convex polytope. Characterizing $\mathcal{C}$ is equivalent to projecting a high-dimensional polytope into the lower-dimensional ${d}$-space, which is computationally impractical due to exponential vertex growth. Other methods such as Fourier-Motzkin Elimination \cite{Khachiyan2009} also suffer from the same combinatorial explosion.

\textit{Baseline method:} We first build a baseline using random sampling. This approach extends a relaxed version of (\ref{eq:040200920}) by making the entries of the degradation vector optimization variables. The discrete optimization variables $x_{i,j}$ are removed and $\bar{\mu}_{i,j}$ become optimization variables.  Additional inequality constraints of the form $0\leq \bar{\mu}_{i,j}\leq (1-d_i)\mu_{i,j}$ reflect the degradation effects. The resulting optimization problem is an LP. Starting from axis-aligned directions, we maximize each $d_i$ individually while keeping others at zero to locate feasible vertices. New points are then randomly generated by averaging feasible pairs and computing the maximum magnitude of the extending rays in their direction under constraints. All such points will belong to the boundary of an approximate baseline capacity region $\hat{\mathcal{C}}_b$.  An untested degradation vector $\hat{{d}}$ is labeled feasible if there exists some ${d} \in \hat{\mathcal{C}}_b$ such that $\hat{{d}} \leq {d}$ (element-wise): \(\bigvee_{{d}\in \hat{\mathcal{C}}_b} \hat{{d}} \leq {d} = \text{true}.\)

\textit{Degradation space exploration:} We make use of the initial step in the baseline approach to collect points on the boundary of $\mathcal{C}$ that maximizes each $d_i$ direction, while all other degradation directions are kept at zero. We augment this set with additional feasible and unfeasible points using the monotonicity rules from equations (\ref{eq:03251431})--(\ref{eq:03251433}) to propagate labels to additional points for free, enabling balancing the classes in the training data.  For a vector $d$ we either add ${d}_{new} = d/2$ (feasible), or ${d}_{new} = d + (1-d)/2$ (infeasible). We use this initial set to train a first classifier. We adopt an entropy-based active learning strategy to sample new points. Let $c({d})$ be the current classifier, represented via the predicted class probability. The next query maximizes uncertainty by solving 
\begin{align}
\label{eq:04091550}
\nonumber
\arg\min_{{d}} &\ c({d})\log c({d}) + (1-c({d}))\log(1-c({d})),
\end{align}
subject to $0\leq {d}\leq 1$. This approach targets points near the boundary. The selected point is tested by querying the feasibility oracle. With these new points and points generated via class balancing, we update the classifier. Uncertainty sampling in active learning~\cite{10.1007/978-1-4471-2099-5_1} includes strategies such as query-by-committee, diversity sampling, or variance reduction~\cite{10.1145/130385.130417,10.5555/1622737.1622744}. While entropy-based sampling may focus too narrowly on borderline cases, this is desirable here, as we aim to delineate the resilience boundary precisely.
\begin{figure}[ht!]
\centering
\includegraphics[width=0.95\linewidth]{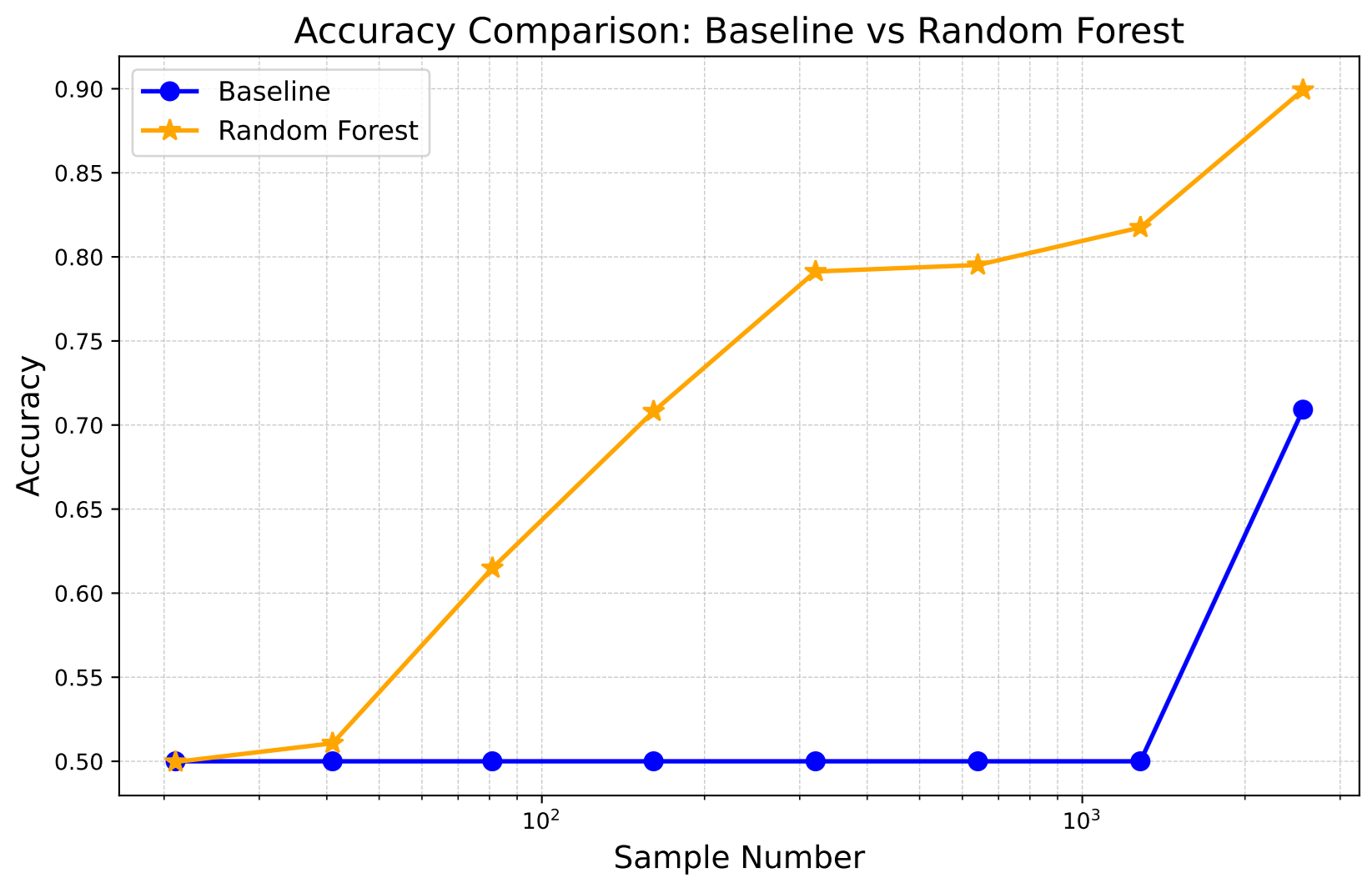}
\caption{Resilience capacity approximation accuracy comparison between the baseline and the active learning approach as a function of the number of training samples used. The results highlight the benefits of active learning (orange), especially when the number of available samples is limited to a few hundred.}
\label{fig:accuracy_comparison}
\end{figure}

\textit{Results:} We compared a random forest (RF) classifier guided by active sampling to the baseline. At each round, entropy is maximized via a gradient-free solver. We train RFs for a set of oracle query budgets ranging from 41 to from 10241 feasibility tests, where samples are generated as described above. RFs were selected over neural networks to reduce architecture-related uncertainty.  We used the {\tt scikit-learn} package \cite{scikit-learn} to fit the RFs with the default settings (e.g., 100 estimators, {\tt gini} criterion for splitting). Fig.~\ref{fig:accuracy_comparison} shows accuracy as a function of the number of training samples. Detailed classification metrics for training budgets of 1281 and 2561 points are reported in Tables~\ref{tab:baseline_scenario1}--\ref{tab:rf_scenario2}. The reports were generated using the {\tt classification\_report} function included in the {\tt scikit-learn} package The reports are based on roughly 10k test points, split equally between feasible and unfeasible classes. The test points were generated by randomly sampling the degradation cube and applying the feasibility test until reaching the target number of test points. Not unexpectedly, the active learning approach achieves superior accuracy with fewer queries.

\begin{table}[ht]
\centering
\caption{Classification Report for Baseline Model (1281 points)}
\label{tab:baseline_scenario1}
\begin{tabular}{lrrrr}
\toprule
 & \textbf{Precision} & \textbf{Recall} & \textbf{F1-Score} & \textbf{Support} \\
\midrule
\textbf{Class 0.0} (unfeasible) & 1.00 & 0.50 & 0.67 & 10834 \\
\textbf{Class 1.0} (feasible) & 0.00 & 0.00 & 0.00 & 0 \\
\midrule
\textbf{Accuracy}    &      &      & 0.50 & 10834 \\
\textbf{Macro avg}   & 0.50 & 0.25 & 0.33 & 10834 \\
\textbf{Weighted avg}& 1.00 & 0.50 & 0.67 & 10834 \\
\bottomrule
\end{tabular}
\end{table}

\begin{table}[ht]
\centering
\caption{Classification Report for Random Forest Model (1281 points)}
\label{tab:rf_scenario1}
\begin{tabular}{lrrrr}
\toprule
 & \textbf{Precision} & \textbf{Recall} & \textbf{F1-Score} & \textbf{Support} \\
\midrule
\textbf{Class 0.0} (unfeasible) & 0.65 & 0.98 & 0.78 & 3602 \\
\textbf{Class 1.0} (feasible) & 0.99 & 0.74 & 0.84 & 7232 \\
\midrule
\textbf{Accuracy}    &      &      & 0.82 & 10834 \\ 
\textbf{Macro avg}   & 0.82 & 0.86 & 0.81 & 10834 \\
\textbf{Weighted avg}& 0.87 & 0.82 & 0.82 & 10834 \\
\bottomrule
\end{tabular}
\end{table}

\begin{table}[ht]
\centering
\caption{Classification Report for Baseline Model (2561 points)}
\label{tab:baseline_scenario2}
\begin{tabular}{lrrrr}
\toprule
 & \textbf{Precision} & \textbf{Recall} & \textbf{F1-Score} & \textbf{Support} \\
\midrule
\textbf{Class 0.0} (unfeasible) & 1.00 & 0.63 & 0.77 & 8568 \\
\textbf{Class 1.0} (feasible) & 0.42 & 1.00 & 0.59 & 2266 \\
\midrule
\textbf{Accuracy}    &      &      & 0.71 & 10834 \\
\textbf{Macro avg}   & 0.71 & 0.82 & 0.68 & 10834 \\
\textbf{Weighted avg}& 0.88 & 0.71 & 0.74 & 10834 \\
\bottomrule
\end{tabular}
\end{table}

\begin{table}[ht]
\centering
\caption{Classification Report for Random Forest Model (2561 points)}
\label{tab:rf_scenario2}
\begin{tabular}{lrrrr}
\toprule
 & \textbf{Precision} & \textbf{Recall} & \textbf{F1-Score} & \textbf{Support} \\
\midrule
\textbf{Class 0.0} (unfeasible) & 0.81 & 0.99 & 0.89 & 4428 \\
\textbf{Class 1.0} (feasible) & 0.99 & 0.84 & 0.91 & 6406 \\
\midrule
\textbf{Accuracy}    &      &      & 0.90 & 10834 \\
\textbf{Macro avg}   & 0.90 & 0.92 & 0.90 & 10834 \\
\textbf{Weighted avg}& 0.92 & 0.90 & 0.91 & 10834 \\
\bottomrule
\end{tabular}
\end{table}

\section{Operational Resilience via Prognostics and Health Management}
\label{sec:phm}

This section demonstrates how the resilience capacity supports operational decisions through PHM steps such as diagnosis, prognostics. We briefly touch on how these steps can be completed by a reconfiguration step. 

\begin{figure}[ht!]
\centering
\includegraphics[width=0.9\linewidth]{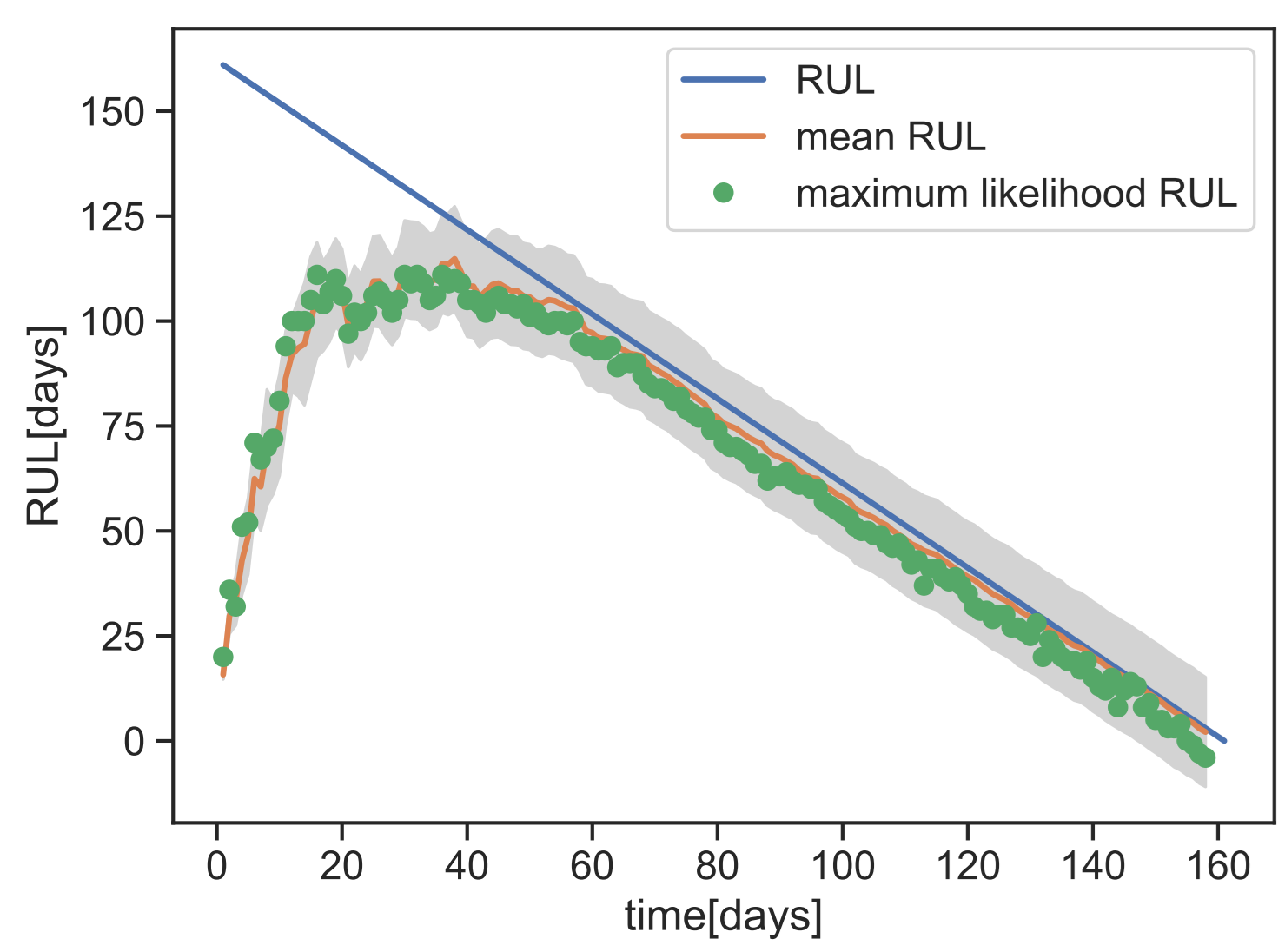}
\caption{Evolution of the system’s RUL over time under machine degradations. Diagnosis results are used to fit the parameters of an exponential degradation model, which is then used to predict future degradation and estimate the RUL. The figure shows the ground truth (blue), mean estimate (orange), maximum‑likelihood estimate (green), and uncertainty range (gray). The uncertainty band reflects both the approximation of the resilience capacity and the parametric degradation model.}
\label{fig:RUL_manufacturing_line}
\vspace{-14pt}
\end{figure}

\textit{Diagnosis} estimates the current system health by identifying degradation magnitudes, which are modeled as reductions in machine capacity. In a model-based setting \cite{DEKLEER200325}, the diagnosis engine uses the system model to detect and isolate faults, returning an estimated degradation vector ${d}$ with associated uncertainties. The resulting degradation estimates are input to the \textit{prognostics} module that forecasts future values of ${d}$ using a degradation model. Models may be physics-based (e.g., Archard wear model~\cite{10.1063/1.1721448}), parametric (e.g., exponential decay~\cite{9263930}), or data-driven (e.g., recurrent neural networks~\cite{10.1145/3631612}). We adopt an exponential model of the form \(d_i(t) = \phi_i + \theta_i e^{\beta_i t + \varepsilon_i - \sigma_i^2/2}\), where $\phi_i$ is constant, $\theta_i$ is log-normally distributed, and $\beta_i$, $\varepsilon_i$ are Gaussian variables. With this type of model we can explicitly compute the likelihood function for fitting the model parameters. Taking logarithms yields a linear-in-parameters model $y_i(t) = \log(d_i(t) - \phi_i) = \hat{\theta}_i + \beta_i t + \varepsilon_i $, with $\hat{\theta}_i = \log \theta_i - \sigma_i^2/2$. Model parameters are estimated using maximum likelihood or MCMC~\cite{Andrieu2003}, based on diagnosis data. Once the forecast model is fit, we evaluate the probability of the (discretized) time when ${d}(t)$ will exit the resilience region. Denoting the time to the boundary of the resilience capacity as $T_{RUL}$, its distribution satisfies: \(    \mathds{P}(T_{RUL} = t_k) = \mathds{E}\left[ \prod_{i=1}^{k-1} c\bigl({d}(t_i)\bigr)\bigl( 1-c\bigl({d}(t_k)\bigr)\bigr) \right]\), where $c(\cdot)$ is the probabilistic classifier for resilience capacity, and the expectation is taken over forecast uncertainty. This estimate includes uncertainty from both the degradation model and the resilience classifier. We simulate exponential degradation for type $B$ machines. The resulting RUL predictions are shown in Fig.~\ref{fig:RUL_manufacturing_line}. Early estimates are uncertain due to limited data but improve over time. We show the mean RUL, maximum likelihood, ground truth, and uncertainty bands based on classifier approximation from 2561 samples.

A system affected by degradations may require \textit{reconfiguration} to adapt its resource allocation and preserve functionality. This can be formulated similarly to the feasibility oracle, but extended over longer horizons to account for future uncertainty in ${d}$. Given long‑term production targets, such a problem can determine if any allocation scheme can satisfy those targets despite anticipated degradations. If delaying maintenance is a priority, the formulation can be adapted to postpone reaching the boundary of the resilience capacity (remaining within $\mathcal{C}$) while still trying to meet production objectives. However, the feasibility oracle must be modified to reflect the stochastic nature of ${d}(t)$. A common approach is to enforce \textit{chance constraints} \cite{Geletu01072013}, ensuring feasibility with high probability. Alternatively, \textit{sampling‑based approximations} \cite{10.5555/2031490} consider multiple scenarios and solve the optimization for each, selecting solutions that minimize expected cost or maximize worst‑case performance. Though accurate, this approach is computationally intensive. A simpler approach replaces ${d}$ with its expected value, yielding a deterministic formulation. For more conservative decisions, robust optimization~\cite{doi:10.1287/opre.1030.0065} finds solutions feasible under worst‑case realizations of ${d}$ within a defined confidence set.

\section{Discussion}
\label{sec:discussion}
The resilience concept proved to be impactful in supply chains and manufacturing, especially in the aftermath of the global disruptions~\cite{OZDEMIR2022101847} induced by the recent COVID-19 pandemic. Although numerous metrics have been proposed to evaluate system resilience, capturing aspects such as output degradation, fault tolerance, and recovery, they remain largely application-specific and lack standardization.

A recent survey of Cyber-Physical Production Systems (CPPS) metrics~\cite{Aruvali2023} highlights a wide range of definitions. These include the ability to maintain production under disruption~\cite{ZHANG2021852}, penalties for system adjustments~\cite{Alexopoulos17122022}, disturbance insensitivity~\cite{Song20}, survival under high-impact events~\cite{CAPUTO2019808}, and adherence to operational thresholds during faults~\cite{su11051447}. These metrics employ diverse tools, e.g., max-plus algebra, fuzzy logic, Monte Carlo simulation, and mixed-integer programming, to evaluate context-specific resilience indicators like throughput or capacity. In contrast, our approach defines resilience with respect to generic, user-defined functional requirements. As long as requirement satisfaction can be evaluated (e.g., through optimization), the method is application-agnostic. That said, domain expertise is still needed to identify degradation mechanisms and construct system models. These models may be built manually or with automated tools, possibly with the help of large language models and generative methods~\cite{maxwell23}. The authors of \cite{9838147} propose a quantitative resilience metric for linear systems under actuator malfunctions, defined in terms of time-to-target, and derive analytical bounds using Lyapunov theory to estimate this metric, emphasizing the key role of redundancy. This aspect is key in our use case as well, and originates in the multi-part functionality of machines. In our approach, the time-to-target would be a requirement, and the uncontrollable actuators would be degradation modes. An interesting perspective on the interplay between robustness and resilience is discussed in \cite{10136414} where resilience is viewed from a safety perspective.

Unlike approaches focusing on isolated failures, our method characterizes the full degradation space and uses active learning to explore it efficiently. We use classifiers to learn the resilience capacity and apply this construct beyond quantification, supporting diagnostic~\cite{Minhas2014}, prognostic~\cite{Bektas_2019}, and reconfiguration planning. Fault detection, diagnostic disambiguation, and predictive maintenance~\cite{Saha2014} also benefit from the same feasibility model. Our prognostics generalize traditional threshold‑based RUL estimates to set‑based boundaries, allowing uncertainty propagation throughout the pipeline. While we use an exponential degradation model, the framework can also accommodate alternatives, such as neural networks~\cite{Bektas2016NARXTS}.

Feasibility testing closely mirrors reconfiguration planning, as both optimize resource allocations under current or anticipated degradations. These tasks align with model predictive control~\cite{GARCIA1989335} and planning strategies~\cite{arshad2005comparison}, sharing a similar constraint structure. The computational cost of the approach stems from three main sources. First, feasibility testing often involves solving large‑scale optimization problems with nonlinear, mixed discrete–continuous dynamics. Surrogate models~\cite{DAVIS2017457} can alleviate some of this load. Second, optimizing large systems is expensive and may benefit from differentiable programming tools~\cite{JIRS_2021}, which enable gradient‑based solvers and parallel execution. Third, learning the resilience capacity involves training classifiers. While decision‑tree methods like random forests are highly interpretable, they scale poorly with data size and dimensionality. In contrast, neural networks better handle high‑dimensional inputs, support GPU parallelism, and enable rapid updating via transfer learning. Together, these components form a flexible but computationally demanding framework for resilience analysis.

\begin{remark}
The main steps of the proposed methodology, applicable across domains, are as follows. The approach starts with a system model and associated simulation engine, along with a set of functional requirements and known degradations or faults. To enable learning of the resilience capacity, the model is augmented for parametric injection of degradations. A simulation budget is then defined, based on the simulation complexity and available training time, to generate training data. An active learning algorithm uses this budget to select informative scenarios, yielding samples for a surrogate model (e.g., RF or NN) that approximates resilience capacity. During operation, this learned model estimates the system’s RUL relative to its requirements. Diagnosis algorithms track degradation states, and the resulting data is used to fit a prognostic model (e.g., exponential decay) that predicts future degradation. Combined with the resilience capacity model, this allows RUL estimation. Finally, reconfiguration algorithms can adjust the system’s configuration or operational settings to extend its RUL. Details on the diagnosis, prognostics, and reconfiguration methods can be found in~\cite{10.1145/3631612}.
\end{remark}

\section{Conclusions}
We presented a methodology for quantifying system resilience by defining and approximating its resilience capacity. Our approach combines classification models and active learning to efficiently characterize high‑dimensional degradation spaces. This framework naturally supports PHM activities such as diagnosis and prognostics, enabling proactive measures to sustain functionality. While applied on a manufacturing case, the approach generalizes across domains. However, several challenges remain. As system complexity grows, differentiable surrogate models can reduce simulation costs, and generative AI and large language models offer promising avenues for automated degradation‑space modeling and requirement extraction. Additionally, many real‑world systems involve hybrid dynamics or nonlinear effects that complicate optimization. Extending this work to incorporate such dynamics, by using specialized solvers or differentiable programming, would further expand its applicability. Moreover, such representations would enable real-time implementations of reconfiguration steps that prolong the RUL under degradations. By connecting resilience quantification with actionable reconfiguration strategies, the proposed framework offers a holistic approach to managing system performance as degradations evolve.

\bibliographystyle{plain}
\bibliography{bibtex/references,bibtex/references_1}  

\end{document}